\documentclass[12pt]{elsarticle}

\usepackage{amsmath}
\usepackage{amssymb}
\usepackage{epsfig}
\usepackage [mathscr]{eucal}
\usepackage{mathrsfs}
\usepackage{xcolor}
\usepackage{booktabs}

\newcommand{\nn}{{\Bbb N}}
\newcommand{\rr}{{\Bbb R}}
\def\beq{\begin{equation}}
\def\eeq{\end{equation}}
\def\be{\begin{eqnarray}}
\def\ee{\end{eqnarray}}
\def\bth{\begin{theorem}}
\def\eth{\end{theorem}}
\def\bde{\begin{definition}}
\def\ede{\end{definition}}
\def\bex{\begin{example}}
\def\eex{\end{example}}
\def\bel{\begin{lemma}}
\def\eel{\end{lemma}}
\def\bep{\begin{proposition}}
\def\eep{\end{proposition}}
\def\bec{\begin{corollary}}
\def\eec{\end{corollary}}

\newtheorem{theorem}{Theorem}[section]
\newtheorem{lemma}{Lemma}[section]
\newtheorem{definition}{Definition}[section]
\newtheorem{corollary}{Corollary}[section]
\newtheorem{proposition}{Proposition}[section]
\newtheorem{remark}{Remark}[section]
\newtheorem{example}{Example}[section]

\makeatletter \@addtoreset{equation}{section} \makeatother

\def\qed{\relax\ifmmode\hskip2em \Box\else\unskip\nobreak\hskip1em $\Box$\fi}

\def\no{\nonumber}
\def\q{\quad}

\begin{document}

\begin{frontmatter}
\title{Domain Recurrence and Probabilistic Analysis of
Residence Time of Stochastic Systems and Domain Aiming Control\tnoteref{label1}}
\author[Rome]{Juliang~Yin\corref{cor1}}
\ead{yin\_juliang@hotmail.com}
\author[Paestum]{Deng~Ding}
\ead{dding@umac.mo}
\author[Khoo]{Suiyang~Khoo}
\ead{sui.khoo@deakin.edu.au}

\cortext[cor1]{Corresponding author. Tel. 86-13662492786}

\address[Rome]{School of Economic and Statistics, Guangzhou University, Guangzhou 510006, China}
\address[Paestum]{Department of Mathematics, University of Macau, Macao}
\address[Khoo]{School of Engineering, Deakin University, VIC 3127, Australia}

\tnotetext[label1]{This work was supported in part by the Natural Science Foundation of China under Grant 61573006.}
\begin{abstract}
The problem of domain aiming control is formulated for controlled stochastic
nonlinear systems. This issue involves regularity of the solution to the resulting closed-loop stochastic system. To begin with, an extended existence and uniqueness theorem for stochastic differential equation with local Lipschitz coefficients is proven by using a
Lyapunov-type function. A Lyapunov-based sufficient condition is also given under which there is no regularity of the solution for a class of stochastic differential equations. The notions of domain recurrence and residence time for stochastic nonlinear systems are introduced, and various criteria for the recurrence and non-recurrence relative to a bounded open domain or an unbounded domain are provided. Furthermore, upper bounds of either the expectation or the moment-generating function of the residence time are derived. In particular, a connection between the mean residence time and a Dirichlet problem is investigated and illustrated with  a numerical example. Finally, the problem of domain aiming control is considered for certain types of nonlinear and linear stochastic systems. Several examples are provided to illustrate the theoretical results.

\end{abstract}
\end{frontmatter}

\section{Introduction}
Let $(\Omega, \mathcal F, \mathbb{P},\{\mathcal F_t\}_{t\ge 0})$ be a complete, filtered probability space, on which a standard $\rr^m$-valued Brownian motion $B_t=(B_t^1, \cdots, B_t^m)^T$ is defined. The filtration $\{\mathcal F_t\}_{t\ge 0}$ is assumed to satisfy the usual conditions, i.e., it is right continuous and $\mathcal F_0$ contains all $\mathbb {P}$-null sets. Let $|x|$ denote the Euclidean norm of $x$ in $\rr^n$, $\rr_+$ the set of nonnegative real numbers and $\nn$ the set of natural numbers. Throughout this paper, if $A $ is a matrix, its transpose is denoted by $A^T $; similarly, the transpose of a vector $x$ is denoted by $x^T$. If $A$ is a symmetric matrix, we denote by
$\lambda_{\max}(A)$ and $\lambda_{\min}(A)$ its largest and smallest eigenvalue, respectively. $a\wedge b$ means the minimum of $a$ and $b$, while $a\vee b$ represents the maximum of $a$ and $b$. Besides, by convention we set $\inf\{\emptyset\}=\infty$, here $\emptyset$ denotes the empty set.

We shall consider a controlled stochastic system described by
\begin{equation}\label{e0}
dX(t)=b(t, X(t), u(t))dt+g(t, X(t), u(t))dB_t, \quad X(0)=x_0, t\ge 0,
\end{equation}
where the state $X(t)\in\rr^n$ and the control input $u(t)\in \rr^{l}$; $b: \rr_+\times \rr^n\times \rr^l\rightarrow\rr^n$
and $g:\rr_+\times \rr^n \times \rr^l\rightarrow\rr^{n\times m}$ are measurable function of $(t, x, u)$. If the control input $u(t)$ is a function of $t$ and $X(t)$, that is $u(t)=u(t, X(t))$, we call it {\sl Markov control}. In particular, if $u(t)$ depends on $X(t)$ only, this is so-called {\sl state feedback control}. A Markov control $u(t)$ is said to be {\sl admissible}, if the solution $\{X^{0, x_0}(t)\}_{t\ge 0}$ of the closed-loop system
\begin{eqnarray}\label{e1}
X(t)&=&x_0+\int_0^t b(s, X(s), u(s, X(s)))ds+\int_0^t g(s, X(s), u(s, X(s)))dB_s\no\\
&\triangleq &x_0+\int_0^t f(s, X(s))ds+\int_0^t \sigma(s, X(s))dB_s, t\ge 0,
\end{eqnarray}
is uniquely determined, that is the equation (\ref{e1}) has a unique global solution in the sense of no explosion. In this case, the unique solution $\{X^{0,x_0}(t)\}_{t\ge 0}$ of the equation (\ref{e1}) is a strong Markov process (see, for example \cite{Bernt:2005}).

Let $U$ be a bounded domain (connected open subset) in $\rr^n$ containing $x=0$ in its interior, and let $\partial U$ be its boundary, $U^c$ the complement and $[U]$ the closure. For any given $T>0, 0<p<1$ and initial value $x_0\in [U]^c$ , if there exists an admissible control $u(t)$ such that the unique solution $X^{0, x_0}(t)$ of the closed-loop system (\ref{e1}) will reach the boundary $\partial U$ during the time interval $[0, T]$ with probability $p$, such a controller is called {\sl domain aiming controller} of $U$. Correspondingly, the controlled stochastic system (\ref{e0}) is also said to be {\sl residence probability controllable} in the domain $[U]^c$.

In the absence of the It\^o integral term in (\ref{e0}), the aforementioned concept of domain aiming control still differs from regional stabilization. The regional stabilization for controlled nonlinear dynamical systems, in general, is to seek an appropriate controller such that the state of the resulting closed-loop system enters a given neighborhood of the origin in a fixed time interval $[0, T]$, but the initial value is confined within a prescribed region of attraction containing the origin \cite {Battilotti:2001}, \cite{Khalil:2002}. In fact, the problem of domain aiming control arises in many controlled stochastic systems, where the goal of control is to accomplish a
 certain task (to ensure that the trajectory of the system reach the target domain) during a specified period ($T$) with a minimum measurement of performance ($0<p<1$). For example, let $Y(t)$ and $Z(t)$ denote the position at time $t$ of a tracker and a tracked object, respectively, then $X(t)=|Y(t)-Z(t)|$ (in fact, $X(t)^2$ is more suitable for mathematical analysis) represents the distance between the tracked one and the tracker. It is very important to find a controller by which the stochastic processes $X(t)$, originating from $x_0$ outside a target domain, goes to the target domain during a period with a probability no less than a given threshold value. As a typical example, in the missile guidance problem \cite{Garnell:1977}, $U$ is defined as the domain of missile interception on a moving target, $T$ is the period of the interception, and $p$ is the probability of successful interception. We would like to stress that the notion of domain aiming control presented in this paper, indeed, also has a significant difference in contrast with the residence time controllability and the aiming control in \cite{Meerkov:1988} and \cite{Kim:1992}. The residence probability controllability in \cite{Kim:1992} gives a stronger formulation
 than the residence time controllability in \cite{Meerkov:1988}.  In particular, the former utilizes the logarithmic first passage probability in a bounded domain to characterize the performance of stochastic systems. In these two references, however, the problem of aiming control and the control technique are suitable for linear systems with small, additive, stochastic perturbations, where the domain is a given bounded set to which the trajectory of the system should be confined during a specified time interval with some probabilistic meanings.

The foundation of domain aiming control is to design an admissible control law such that the resulting closed-loop system has a regional stability in some probabilistic sense. In order to describe the behavior of the solution to (\ref{e1}) originating from outside of a domain $U$, the first hitting time of the domain plays an important role, which can be seen as the residence time outside the domain $U$. However, since the explicit solution to (\ref{e1}) is generally unavailable, it is often difficult or even impossible to get the probability distribution of the residence time. Alternatively, Lyapunov-type conditions are given for obtaining an upper bound to the probability that the residence time is not greater than the period $T$. This idea is not new. In fact, Lyapunov-type methods have been used to obtain upper bounds to probabilities of certain events for dealing with stabilization in probability of nonlinear stochastic systems in \cite{Husher:1966} and \cite{Battilotti:2003}.

As stated above, the domain aiming control problem involves the regularity of the resulting closed-loop stochastic system. A well-known uniqueness and existence theorem of solutions for a stochastic differential equation driven by a Brownian motion requires that the coefficients of the equation satisfy the local Lipschitz condition and the linear growth condition. The linear growth condition is too restrictive to be satisfied in many cases for domain aiming control problems of stochastic systems. In view of this limitation, an extended existence and uniqueness theorem for stochastic differential equations is proven by using a
Lyapunov-type function analogous to Khasminskii's criteria in \cite{Hasminskii:2012}. Lyapunov-based sufficient conditions are also given, under which there is no regularity of the solution for  a class of stochastic differential equations. For a regular stochastic system, recurrence is an extremely important concept for studying domain aiming control problems. Roughly speaking, the recurrence relative to a domain means that the trajectory issuing from any initial values outside the domain will eventually reach the domain with probability one. The mathematical definition of the recurrence will be given in Section 3.
If the solution of a stochastic differential equation is recurrent relative to a domain $U$, then the residence time of the solution outside the domain is almost surely finite. In order to deal with domain aiming control problems, we need to compute the probability that the residence time is less than or equal to the duration period.
If the target domain is not recurrent, it is impossible to discuss such a control problem provided the probability $p$ is relatively large. Therefore, sufficient conditions for domain recurrence and the finiteness of the expected residence time are given by virtue of Lyapunov functions. In addition, sufficient conditions are also proposed for no recurrence of stochastic systems. In domain aiming control problems, the probability that the trajectory $\{X^{0,x_0}(t)\}$ first reaches the boundary $\partial U$ during the time interval $[0, T]$ is needed to estimate, because a direct calculation for this is usually infeasible. We derive a lower bound  of the probability by obtaining upper bounds of its inverse probability or the mean residence time, and Chebyshev's inequality is used to get the required estimates. Under
suitable assumptions, it is shown that there exists a connection between the mean residence time and a Dirichlet problem. Although the partial differential equation in the Dirichlet problem, in general, is not explicitly solvable, numerical methods can be used to give an approximate result of the solution.
An illustrative example is provided for this purpose. By using these fundamental results, some sufficient conditions are given under which a class of nonlinear  controlled stochastic systems is domain aiming controllable by selecting a suitable admissible control law. For the case of linear controlled stochastic systems,
the domain aiming controllability follows by designing a feedback controller under weaker conditions.

If a nonlinear dynamical system incorporates a linear nominal part perturbed by model uncertainties, nonlinearities and both additive and multiplicative random noise, modeled by a Brownian motion, an It\^o stochastic differential equation is suitable for describing such a real system in engineer. The basic theory and applications of It\^o stochastic integrals and stochastic differential equations have been given in many textbooks (see, for example \cite{Friedman:1975}, \cite{Karatzas:1999}, \cite{Bernt:2005} and so on ), while various types of stochastic stability for stochastic systems
have been discussed in \cite{Husher:1967}, \cite{Hasminskii:2012} and \cite{Mao:2008}. In particular, a type of finite-time stochastic stability for stochastic processes were given in \cite{Husher:1967}, which is associated with certain types of first exit time problems. This stability means that a stochastic process $x_t$ will remain within in a given region $Q_2$ in an interval $[0, T]$ with
probability no less than $1-\lambda$, if $x_0=x\in Q_1\subset Q_2$. In this paper, the residence probability controllability of a controlled stochastic system corresponds to that the closed-loop system admits a regional stability in probability. This stochastic stability, in connection with the first hitting time problem, can be defined in a similar way as finite-time stochastic stability in \cite{Husher:1967}. Obviously, this kind of stochastic stability is different from almost surely finite-time
stability of stochastic systems in \cite {Yin:2011} and \cite{Khoo2013}.

\par
The paper is organized as follows. Section 2 begins with the regularity of the equation (\ref{e1}) with the help of Lyapunov functions. A sufficient condition is also given, under which there is no regularity for the equation. Section 3 devotes to giving criteria for the recurrence and non-recurrence of the solution to the equation (\ref{e1}) relative to a bounded open domain. Upper bounds of either the expectation or the moment-generating function of the residence time will be derived. In particular, a connection between the mean residence time and a Dirichlet problem has been built. A few examples are provided to illustrate the theoretical results including a numerical example. In section 4, the problem of domain aiming control is considered for nonlinear stochastic systems and linear stochastic systems. Finally, some concluding remarks are given in Section 5.

\section{Regularity of the solution}

\par
Let us first consider the existence and uniqueness of the solution to (\ref{e1}) under unrestrictive conditions. A typical existence and uniqueness theorem states that there exists a unique solution to the equation (\ref{e1}) if both $f$ and $\sigma$ satisfy the local Lipschitz condition and the linear growth condition \cite {
Friedman:1975}, \cite{Hasminskii:2012}, \cite{Mao:2008}. The linear growth condition is rather restrictive, so we shall impose some mild conditions on the coefficients of the equation. It is clear that the equation admits a unique local solution $\{X^{0,x_0}(t)\}$ for any $x_0\in \rr^n$ if the local Lipschitz condition is fulfilled. In this case, we shall say that the process $\{X^{0,x_0}(t)\}$ is {\sl regular} provided the explosion time is almost surely infinite. The regularity implies that the equation has a unique global solution. In the sequel, unless explicitly stated, we will denote $ X^{0,x_0}(t)$ by $X(t)$ without stressing the initial value $x_0$ for simplicity's sake.

\par
Let $\mathcal L^1(\rr_+, \rr_+)$ denote the set of all functions $\nu: \rr_+\rightarrow \rr_+$ such that $\int_0^\infty \nu(t)dt<\infty$. We also denote by $\mathbf C^{1,2}$ the set of all functions $V$ on $\rr_+\times \rr^n$ which are once differentiable in $t\in\rr_+$ and continuously twice differentiable in $x\in\rr^n$. If $V\in \mathbf C^{1,2}$, we define an operator $\mathcal L $ acting on $V$ by
\begin{equation*}
\mathcal LV(t, x)=\frac{\partial V(t,x)}{\partial t}+\sum_{i=1}^n f_i(t, x)\frac{\partial V(t,x)}{\partial x_i}+\frac{1}{2}\sum_{i,j=1}^n a_{ij}(t,x)\frac{\partial^2 V(t,x)}{\partial x_i\partial x_j},
\end{equation*}
where $a_{ij}(t, x)=(\sigma(t,x)\sigma(t,x)^T)_{i, j}$.
\par

\begin{theorem}\label{l1} Suppose that $f$ and $\sigma$ are all continuous in $(t,x)\in \rr_+\times \rr^n $ and satisfy the following local Lipschitz condition: For each integer $n\ge 1$ and arbitrary $T>0$, there exists a positive constant $K_{n, T}$ such that, for all $t\in [0, T]$, as $
|x|\vee|y|\leq n$,

\begin{equation}\label{e2}
|f(t,x)-f(t,y)|\leq K_{n, T}|x-y|, \q |\sigma(t, x)-\sigma(t,y)|\leq K_{n, T}|x-y|.
\end{equation}
If there exists a nonnegative function $V\in \mathbf C^{1, 2}$ such that
\begin{eqnarray}
&&\mathcal LV(t,x)\leq \gamma(t)+\alpha(t) V(t, x), \q \forall (t,x)\in [0, T]\times\rr^n, \label{e3}\\
&& \inf_{0\leq t\leq T, |x|=R} V(t, x)\rightarrow \infty, R\rightarrow\infty,\label{e3-1}
\end{eqnarray}
where $\gamma(\cdot)$ and $\alpha(\cdot)$ are nonnegative functions satisfying
$$
\int_0^T \left(\gamma(t)+\alpha(t)\right) dt<\infty, \forall T>0,
$$
then the equation {\rm(\ref{e1})} has a unique global solution.
\end{theorem}
\par
{\bf Proof.} For each $n\ge 1$, let
$$
f_n(t, x)=\left\{ \begin{array}{lll} f(t,x)\q\q &\mbox{if}\q & |x|\leq n,\\
f(t, nx/|x|), \q\q &\mbox{if}\q &|x|>n,
\end{array}
\right.
$$
and $\sigma_n(t,x)$ is similarly defined. Since $f(t, 0)$ and $\sigma(t, 0)$ are continuous functions in $t\in [0, T]$, it is readily seen that $f_n$ and $\sigma_n$ satisfy the Lipschitz condition and the linear growth condition. Hence there is a unique solution $X_n(t)\in \mathcal M^2([0, T]; \rr ^n)$ solving the following stochastic differential equation
$$
dX_n(t)=f_n(t, X_n(t))dt+\sigma_n(t, X_n(t))dB_t, \q X_n(0)=x_0, \q t\in [0, T],
$$
where $\mathcal M^2([0, T]; \rr ^n)$ denotes the class of all adapted processes $f(t)$ satisfying $\mathbb {E}\int_0^T |f(t)|^2 dt<\infty$.
Let $\tau_n$ be the first exit time of the processes $X_n(t)$ from the set $\{x\in \rr^n; |x|<n\}$ in $[0, T]$. It is clear that $\tau_n$ is a stopping time, which is defined as
$$
\tau_n=\inf\{t\in[0, T]: |X_n(t)|\ge n\}\wedge T.
$$
It can be shown (see Friedman \cite{Friedman:1975}, Theorem 2.1, p. 102]) that
\begin{equation}\label{e4}
X_n(t)=X_{n+1}(t) \;\; \mbox{a.s.} \q \mbox{if}\q t\in [0, \tau_n].
\end{equation}
This implies that $\tau_n$ is increasing and has a limit $\tau=\lim_{n\rightarrow\infty}\tau_n$. Let $X(t)=X_n(t)$ as $t\in[0, \tau_n]$. Then by the uniqueness result (\ref{e4}), the processes $X(t)$ is well defined as $t\in[0, \tau_n]$. Next, we will show that $\mathbb {P}(\tau=T)=1$. If this is not true, then there exists a sufficiently small positive constant $\epsilon$ such that $\mathbb {P}(\tau\leq T-\epsilon)\ge 2\epsilon$. Hence, there exists a sufficiently large integer $n_0$, as $n\ge n_0$, $\mathbb {P}(\tau_n\leq T-\epsilon)\ge \epsilon$. In this situation, for $0\leq t\leq T$, we can apply It\^o's formula to derive that
\begin{eqnarray*}
& &\mathbb{ E} V(t\wedge \tau_n, X(t\wedge\tau_n))=V(0, x_0)+\mathbb {E}\int_0^{t\wedge\tau_n} \mathcal LV(s, X(s))ds\\
& &\q \leq V(0, x_0)+\int_0^T \gamma(s)ds+\mathbb {E}\int_0^{t\wedge \tau_n}\alpha(s)V(s, X(s))ds\\
& &\q \leq V(0, x_0)+\int_0^T \gamma(s)ds+\int_0^{t}\alpha(s)\mathbb {E} V(s\wedge \tau_n, X(s\wedge\tau_n))ds.
\end{eqnarray*}
It then follows from Gronwall's inequality that
$$
\mathbb {E} V(t\wedge \tau_n, X(t\wedge\tau_n))\leq \left(V(0, x_0)+\int_0^T \gamma(s)ds \right)\exp\left(\int_0^T \alpha(s)ds \right)<\infty.
$$
Now fixing $t=T-\epsilon$, we then obtain
\begin{eqnarray*}
\epsilon \inf_{0\leq t\leq T, |x|=n}V(t, x)&\leq& \mathbb {P}(\tau_n\leq T-\epsilon) \inf_{0\leq t\leq T-\epsilon, |x|=n}V(t, x)\\
&\leq &\left(V(0, x_0)+\int_0^T \gamma(s)ds \right)\exp\left(\int_0^T \alpha(s)ds \right)
\end{eqnarray*}
This is a contradiction since $\inf_{0\leq t\leq T, |x|=n}V(t, x)\rightarrow\infty$ as $n\rightarrow\infty$. We thus have that
$\mathbb {P}(\tau=T)=1$. Since $T$ is an arbitrary positive constant, $\{X(t)\}_{t\ge 0}$ is a unique global solution of (\ref{e1}).\hfill $\Box$
\par
From the proof of Theorem \ref{l1}, it is known that the processes $X(t)$ belongs to $\mathcal M^2([0, T]; \rr^n)$ for all $t\in [0, T], \forall T>0$. In addition, if condition (\ref{e3}) is replaced by the following
$$
\mathcal LV(t,x)\leq \gamma(t)+\alpha(t)\left(V(t, x)\right)^{r}, \q r\in[0,1),
$$
and the other assumptions remain unchanged, then Theorem \ref{l1} still holds true. Indeed, by using a Bihari-type inequality (see Mao \cite{Mao:2008}, Theorem 8.3, p. 46), we can get that
\begin{eqnarray*}
&&\mathbb {E}\big( V(t\wedge\tau_n, X(t\wedge\tau_n))\big)\\
&&\leq \left(\Big(V(0, x_0)+\int_0^T \gamma(s)ds \Big)^{1-r}+(1-r)\int_0^T\alpha(s)ds\right)^{\frac{1}{1-r}}<\infty.
\end{eqnarray*}
By this inequality, the remainder of the proof follows in the same way as that of Theorem \ref{l1}.
\begin{remark}
We remark that Theorem \ref{l1} has generalized some existing results in the literature. To be precise, if $\gamma(t)=0$ and $\alpha(t)=\alpha>0$, then $\mathcal LV(t, x)\leq \alpha V(t,x)$, which has been considered by Hasminskii {\rm {(see \cite {Hasminskii:2012}, Theorem 3.5, p.75)}}. Assume that $f$ and $\sigma$, except the local Lipschitz condition {\rm{(\ref{e2})}}, also satisfies the following monotone condition: For any $T>0$, there exists a positive constant $K_T>0$ such that for all $(t, x)\in [0, T]\times\rr^n$,
\begin{equation}
x^T f(t,x)+\frac{1}{2}|\sigma(t,x)|^2\leq K_T(1+|x|^2).
\end{equation}
By taking $V(t, x)=|x|^2$, one can easily
derive that
$$
\mathcal LV(t, x)\leq 2K_T+2K_T V(t,x).
$$
This corresponds to the case of $\gamma(t)=\alpha(t)\equiv 2K_T$. Hence the existence and uniqueness result of the global solution in \cite{Mao:2008} (see Theorem 3.6, p. 59) can be seen as a special case of Theorem \ref{l1}.
\end{remark}
\begin{remark} \label{r1} Suppose that there exists an $R_0>0$, for any $T>0$ the condition {\rm(\ref{e3})} is weakened as
$$
\mathcal LV(t,x)\leq \gamma(t)+\alpha(t) V(t, x), \q \forall (t,x)\in [0, T]\times\{x\in\rr^n: |x|\ge R_0\},
$$
and other assumptions remain unchanged, then Theorem {\rm \ref{l1}} is still valid. Without loss of generality, assume that $x_0\in U_{R_0}=\{x\in \rr^n: |x|<R_0\} $. We now choose a sufficiently large $n\in \mathbb{N}$ such that $U_{R_0}\subset \{x\in\rr^n: |x|<n\}$. From the proof of Theorem \ref{l1}, it is known that
$X(t)=X_n(t)$ is the unique solution of {\rm(\ref{e1})} as $t\in [0, \tau_n]$. Let $\tau_{R_0}$ is the first passage time of the process $X(t)$ from the domain $ U_{R_0}$
in $[0, T]$. If $\mathbb {P}(\tau_{R_0}=T)=1$, then $\mathbb {P}(\tau_n=T)=1$ due to $\tau_{R_0}\leq \tau_n$ a.s. This means that $\mathbb P(\tau=T)=1$. If $\mathbb{P}(\tau_{R_0}=T)<1$, by It\^o's formula and Gronwall's inequality, we can derive
\begin{eqnarray*}
\mathbb {E} V(t\wedge \tau_n, X(t\wedge\tau_n))&\leq &\left(\mathbb{E}V(t\wedge\tau_{R_0}, X(t\wedge\tau_{R_0}))+\int_0^T \gamma(s)ds \right)\\
& &\q \times \exp\left(\int_0^T \alpha(s)ds \right).
\end{eqnarray*}
It is easily seen that
\begin{equation*}
\mathbb {E}(V(t\wedge\tau_{R_0}, X(t\wedge\tau_{R_0}))\leq \sup_{0\leq t\leq T, |x|\leq R_0} V(t, x)<\infty.
\end{equation*}
The remained proof is similar, so we omit the details.
\end{remark}
\par

The following theorem gives a sufficient condition under which the regularity of the solution to the equation (\ref{e2}) does not hold. This theorem also generalizes a result of Khasminskii \cite{Hasminskii:2012}, and has an analogous proof.
\begin{theorem} Under the condition of {\rm{(\ref{e2})}}, if there exists a nonnegative and bounded function $V\in \mathbf C^{1,2}$ in $\rr_+\times \rr^n$ such that
\begin{equation}\label{e5}
\mathcal LV(t,x)\ge \bar\alpha(t) V(t,x), \q V(0,x_0)>0,
\end{equation}
where $\bar\alpha(t)\ge 0$ and satisfies
$$
\int_0^T \bar \alpha(t)dt<\infty, \forall T>0, \q \lim_{T\rightarrow\infty}\int_0^T \bar \alpha(t)dt\rightarrow\infty,
$$
then the solution of the equation {\rm (\ref{e2})} is not regular.
\end{theorem}
\par
{\bf Proof.} As in the proof of Theorem 2.1, define $X(t)=X_n(t)$ as $t\in [0, \tau_n]$. We will prove that $\mathbb {P}(\tau=T)<1$, where $\tau=\lim_{n\rightarrow\infty}\tau_n$. For $0\leq t\leq T$, we apply It\^o's formula to $e^{-\int_0^t \bar \alpha(s)ds} V(t, X(t))$ and utilize (\ref{e5}) to obtain $$
\mathbb {E}\left(e^{-\int_0^{t\wedge \tau_n}\bar\alpha(s)ds} V\big(t\wedge \tau_n, X(t\wedge \tau_n)\big)\right)\ge V(0,x_0).
$$
Since $V$ is bounded, it follows from the Lebesgue dominated convergence theorem that
\begin{equation}\label{e6}
\mathbb {E}\left(e^{-\int_0^{t\wedge \tau}\bar\alpha(s)ds} V\big(t\wedge \tau, X(t\wedge \tau)\big)\right)\ge V(0,x_0).
\end{equation}
Assume that $\mathbb {P}(\tau=T)=1$, we will derive a contradiction. Indeed, by (\ref{e6}), one may easily obtain that
$$
e^{-\int_0^t \bar \alpha(s)ds}\sup_{\rr_+\times \rr^n} V(t, x)\ge V(0, x_0), \forall t\in [0, T],
$$
which implies that
$$
\int_0^t \bar \alpha(s)ds\leq \log \frac{\sup_{\rr_+\times \rr^n} V(t, x)}{V(0, x_0)}<\infty, \forall t\in [0, T].
$$
This is in contradiction with the properties of $\lim_{t\rightarrow\infty}\int_0^t\bar\alpha(s)ds\rightarrow\infty$ and the boundedness of $V$, since $T$ is an arbitrary positive constant.
This completes the proof. \hfill $\Box$

\section{Domain recurrence and residence time of stochastic systems}
 The notion of recurrence relative to a domain was introduced by Khasminskii \cite{Hasminskii:2012}. The domain recurrence of a stochastic system implies that its trajectories originating from any initial values outside the domain will finally enter the domain with probability one. The following domain recurrence requires almost all the paths of the system will eventually reach the boundary of the domain. There is a slight difference between the two definitions (Definition 3.1 corresponds to the recurrence of $[D]$ indeed). We also provide Lyapunov-type sufficient conditions for a regular stochastic system to be recurrent or non-recurrent relative to an open bounded domain. In addition, we derive upper bounds of the mean and moment-generating function of the residence time under certain conditions, by which the probability that the trajectories first reach the boundary of the domain in fixed time achieves an upper bound. Under suitable assumptions, we will give a connection between the expectation of the stochastic residence time and a Dirichlet problem.

\begin{definition} Let $D$ be a (bounded or unbounded) domain. A Markov process $X(t)$ is said to be recurrent relative to the domain $D$ (or $D$-recurrent) if it is regular, and for every $s\in \rr_+, x\in D^c$,
$$
\mathbb {P}^{s, x}(\tau_D<\infty)=1,
$$
where $\tau_D=\inf\{t\ge s: X(t)\in \partial D, X(0)=x\}$, called the residence time in the domain $[D]^c$, is the first hitting time of $\partial D$, and $\mathbb {P}^{s, x}$ denotes the probability law of $\{X(t)\}_{t\ge s}$ when its initial value is $X(s)=x$. If $s=0$, for simplicity, we write $\mathbb {P}^{s,x}$ as $\mathbb {P}^{x}$.
\end{definition}

\begin{definition}
A function $\mu: \rr_+\rightarrow \rr_+$ is said to belong to class $\mathcal K$ if it is
continuous, strictly increasing and satisfies $\mu(0)=0$. A class $\mathcal K$
function $\mu$ is said to belong to class $\mathcal K_\infty$ if
$\mu(r)\rightarrow \infty$ as $r\rightarrow\infty.$

\end{definition}

\begin{theorem} \label {T1}Let $U$ be an open bounded domain with $0$ in its interior. Suppose that $f$ and $\sigma$ are all continuous in $(t,x)\in \rr_+\times \rr^n $ and satisfy {\rm (\ref{e2})}. If there exists a nonnegative function $V\in \mathbf C^{1,2}$ such that,
\begin{equation}\label{e6-1}
\inf_{t\ge 0, |x|\ge R} V(t, x)\rightarrow \infty, \q R\rightarrow\infty,
\end{equation}
and
\begin{eqnarray}\label{e6-2}
\mathcal LV(t,x)\leq \nu(t)\wedge\left(\nu(t)+|\frac{\partial V}{\partial x}(t, x)\sigma(t,x)|^2-\mu(|x|)\right),
\end{eqnarray}
where $\nu\in \mathcal L^1(\rr_+, \rr_+)$ and $\mu\in \mathcal K$, then the regularity holds for the equation {\rm(\ref{e1})} and its solution $X(t)$  is recurrent with respect to the domain $U$ for any initial value $x_0\not\in U$. Furthermore, if {\rm (\ref{e6-2})} is strengthened as
\begin{equation} \label{e6-3}
\mathcal LV(t,x)\leq \nu(t)-\mu(|x|),
\end{equation}
and assume that $\tau_U$ is the residence time of the sample path for the equation {\rm(\ref{e1})} outside the domain $U$ (in the domain $[U]^c$), then an upper bound of $\mathbb {E}^{x_0}(\tau_U)$ is given by
\begin{equation} \label{e6-4}
\mathbb {E}^{x_0}(\tau_U)\leq \frac{V(0, x_0)+\int_0^\infty \nu(t)dt-\inf_{t\ge 0, x\in \partial U}V(t,x)}{\mu(r)},
\end{equation}
where $r=\inf_{x\in\partial U}|x|$ and $\mathbb {E}^{x_0}$ denotes the expectation with respect to the probability law $\mathbb {P}^{x_0}$.
\end{theorem}
\par
{\bf Proof.} The regularity follows directly from Theorem \ref{l1} with $\gamma(t)=\nu(t)$ and $\alpha(t)=0$. It remains to show that the solution $X(t)$ is $U$-recurrent for any given initial value $x_0\not\in [U]$, since the conclusion is trivial in the case of $x_0\in \partial U$. Let $r=\inf_{x\in\partial U}|x|$. It is clear that $0<r<\infty$. For any $\delta>0$, define $U_\delta=\{x\in \rr^n: |x|<\delta\}$. We now choose a sufficiently large $R>0$ so that $U\subset U_R$ and $x_0\in U_R$. Define the stopping time
$$
\tau_R=\inf\{ t\ge 0: |X(t)|\ge R\}=\inf\{t\ge 0: X(t)\in U_R^c\}.
$$
An application of It\^o's formula yields that
\begin{eqnarray*}
V(t, X(t))&=&V(0, x_0)+\int_0^t \nu(s)ds-\int_0^t \big(\nu(s)-\mathcal LV(s, X(s))\big)ds\\
& &+\int_0^t \frac{\partial V}{\partial x}(s, X(s))\sigma(s, X(s))dB_s.
\end{eqnarray*}
Let
$$
M(t)=\int_0^t \frac{\partial V}{\partial x}(s, X(s))\sigma(s, X(s))dB_s.
$$
By using (\ref{e6-2}) and an convergence theorem of semi-martingales established by Liptser and Shiryayev \cite{Liptser:1989}, we have $\lim_{t\rightarrow \infty} V(t, X(t))$ exists and is finite almost surely and, moreover,
$$
\int_0^\infty \big (\nu(t)-\mathcal LV(t, X(t))\big)dt<\infty, \q \mbox{a.s.}
$$
and
$$
-\infty<\lim_{t\rightarrow \infty} M(t)<\infty, \q \mbox{a.s.,}
$$
due to $\nu\in \mathcal L^1(\rr_+, \rr_+)$ and $\nu(s)-\mathcal LV(s,x)\ge 0$. Furthermore, it can be shown by adopting the same approach in
\cite{Mao:2001} that,
$$
\int_0^\infty \big|\frac{\partial V}{\partial x}(t, X(t))\sigma(t, X(t))\big|^2 dt<\infty \q \mbox{a.s.,}
$$
which also implies that
\begin{eqnarray*}
&&\int_0^\infty \mu(|X(t)|)dt\\
&&\leq \int_0^\infty \big(\nu(t)-\mathcal LV(t, X(t))+\big|\frac{\partial V}{\partial x}(t, X(t))\sigma(t, X(t))\big|^2\big) dt<\infty
\end{eqnarray*}
holds with probability 1. By using stochastic Barb\u{a}lat's lemma in \cite{Yu:2010} (see Lemma 3, where the assertion can be proved in the same way as that in \cite{Deng:2001}), we have
$$
\mathbb {P}^{x_0}(\lim_{t\rightarrow\infty}|X(t)|=0)=1,
$$
since $\mu(\cdot)$ is continuous and positive definite. Therefore, the recurrence of the solution to (\ref{e1}) relative to the domain $U$ follows. In order to give an upper bound for $\mathbb E^{x_0}(\tau_U)$, we now use It\^o's formula and (\ref{e6-3}) to get
\begin{eqnarray}
&&\mathbb {E}^{x_0} V(\tau_U\wedge\tau_R\wedge t, X(\tau_U\wedge\tau_R\wedge t))\no\\
 & &\q\leq V(0, x_0)+\int_0^t \nu(s)ds-\mathbb {E}^{x_0}\int_0^{\tau_U\wedge\tau_R\wedge t}\mu(|X(s)|)ds\nonumber\\
 &&\q \leq V(0, x_0)+\int_0^t \nu(s)ds-\mu(r) \mathbb {E}^{x_0}\big(\tau_U\wedge\tau_R\wedge t\big).\label{e7}
\end{eqnarray}
Letting $R\rightarrow\infty$ and applying Fatou's lemma on both sides of (\ref{e7}), the regularity means that
$$
\mathbb {E}^{x_0} V(\tau_U\wedge t, X(\tau_U\wedge t))\leq V(0, x_0)+\int_0^\infty \nu(t)dt-\mu(r) \mathbb {E}^{x_0}\big(\tau_U\wedge t\big).
$$
Since $V$ is a nonnegative function, it is easy to derive that
$$
\mathbb {E}^{x_0}(\tau_U)\leq \frac{ V(0, x_0)+\int_0^\infty \nu(t)dt}{\mu(r)}<\infty.
$$
An application of Fatou's lemma again, together with the monotone convergence theorem, yields that
$$
\mathbb {E}^{x_0} V(\tau_U, X(\tau_U))\leq V(0, x_0)+\int_0^\infty \nu(t)dt-\mu(r) \mathbb {E}^{x_0}\big(\tau_U\big)
$$
This means that
$$
\mathbb {E}^{x_0}\big(\tau_U\big)\leq \frac{V(0, x_0)+\int_0^\infty \nu(t)dt-\inf_{t\ge 0, x\in \partial U}V(t,x) }{\mu(r)}.
$$
The proof is complete. \hfill $\Box$

\begin{remark}\label{R2}
It is known from the proof of Theorem \ref{T1} that if we can find a nonnegative function $V\in \mathbf C^{1,2}$ satisfying {\rm(\ref{e6-1})} and {\rm(\ref{e6-3})} in the domain $\rr_+\times U^c$, then the solution of {\rm(\ref{e1})} is $U$-recurrent. Indeed, the regularity holds according to Remark {\ref{r1}}. Although the condition {\rm (\ref{e6-3})} is required to be satisfied only in the domain $\rr_+\times U^c$, it not only can ensures the recurrence of the domain $U$ but also can derive an upper bound for $\mathbb {E}^{x_0}(\tau_U)$. It is obvious that the condition {\rm (\ref{e6-3})} is restrictive for the recurrence of the domain $U$. To relax this condition, we will consider a special case: For any $R>0$, if $U\subset U_R$ and $x_0\in U_R\setminus U$, then the condition {\rm (\ref{e6-3})} can be weakened as
$$
\mathcal LV(t, x)\leq \nu(t), \q (t,x)\in \rr_+\times U^c,
$$
provided that the solution of {\rm(\ref{e1})} could leave the domain $U_R\setminus U$ in finite time with probability one, namely, $\mathbb P^{x_0}(\tau_U\wedge\tau_R<\infty)=1$. In this case, by It\^o's formula, we have that
\begin{eqnarray*}
\mathbb{E}^{x_0}\Big(V(t\wedge\tau_U\wedge\tau_R, X(t\wedge\tau_U\wedge\tau_R))\Big)
\leq V(0, x_0)+\int_0^{t}\nu(s)ds.
\end{eqnarray*}
If then follows from the assumption of $\mathbb P^{x_0}(\tau_U\wedge\tau_R<\infty)=1$ and Fatou's lemma that

$$
\mathbb{E}^{x_0} \Big(V(\tau_U\wedge\tau_R, X(\tau_U\wedge\tau_R))\Big)\leq V(0, x_0)+\int_0^{\infty}\nu(t)dt\triangleq K<\infty,
$$
which also gives that
$$
\mathbb{P}^{x_0}(\tau_R<\tau_U)\inf_{t\ge 0, |x|\ge R} V(t, x)\leq \mathbb{E}^{x_0} \Big(V(\tau_U\wedge\tau_R, X(\tau_U\wedge\tau_R))\Big)\leq K.
$$
Therefore, by taking $R\rightarrow \infty$ and using the regularity of solution and {\rm(\ref{e6-1})},
$$
\mathbb{P}^{x_0}(\tau_U<\infty)=1.
$$
However, the assumption of $\mathbb P^{x_0}(\tau_U\wedge\tau_R<\infty)=1$ is rather rarely satisfied. As we shall see in Theorem \ref{T2}, conditions of {\rm{(\ref{e11})}} or {\rm{(\ref{e12})}} are sufficient for this assumption to hold true.
\end{remark}

It is worth pointing out that, if there exist two $\mathcal K$ functions $\mu_1$ and $\mu_2$ such that
$\mu_1(|x|)\leq V(t, x)\leq \mu_2(|x|)$, and (\ref{e6-2}) is satisfied, then the stochastic LaSalle theorem in \cite{Mao:2001} yields that
$\mathbb {P}^{x_0}(\lim_{t\rightarrow \infty}X(t, x_0)=0)$, and the solution of (\ref{e1}) is recurrent with respect to the domain $U$.
Therefore, the condition (\ref{e6-1}) in Theorem \ref{T1} is slightly mild for the recurrence of a bounded domain.

Theorem \ref{l1} and Theorem \ref{T1} depend on different auxiliary functions, so it is very necessary to point out how to construct auxiliary functions $V_1$ and $V_2$ that satisfy (\ref{e3}), (\ref{e3-1}) and (\ref{e6-1}), (\ref{e6-2}) for stochastic differential equations, respectively. We will give an example to illustrate these issues.

\begin{example} Let $U=(-\delta, \delta), \delta>0$. Let us consider the following one-dimensional stochastic system
\begin{equation} \label{e8}
dX(t)=f(X(t))dt+\sigma(X(t))dB_t, \q t\ge 0,
\end{equation}
with initial value $X(0)=x_0\not\in U$. Assume that
$$
f(x)=\alpha_1 x, \q \sigma(x)=x\sqrt{\alpha_2+\alpha_3 x^2}, \q \alpha_2, \alpha_3\ge 0.
$$
It is obvious that $f$ and $\sigma$ satisfy the local Lipschitz condition. Let $V_1(x)=\log(1+x^2)$. By It\^o's formula, we have
$$
\mathcal LV_1(x)=\frac{2\alpha_1 x^2}{1+x^2}+\frac{x^2(\alpha_2+\alpha_3 x^2)}{1+x^2}-\frac{2x^4(\alpha_2+\alpha_3x^2)}{(1+x^2)^2}\leq 2|\alpha_1|+\alpha_2\vee \alpha_3.
$$
For $x\in U^c=\{x\in \rr: |x|\ge \delta\}$, define $V_2(x)=K+\log|x|$, where $K$ is a sufficiently large positive constant such that
$V_2(x)\ge 0$. It follows from It\^o's formula that
$$
\mathcal LV_2(x)=\frac{f(x)}{x}-\frac{1}{2}\frac{\sigma^2(x)}{x^2}=\alpha_1-\frac{\alpha_2}{2}-\frac{\alpha_3}{2} x^2.
$$
If the constants $\alpha_1, \alpha_2$ and $\alpha_3$ have the property of
$$
\frac{\alpha_3^2}{2}\delta^2>\alpha_1-\frac{\alpha_2}{2}, \alpha_2\ge 0, \alpha_3>0,
$$
then by Theorem \ref{l1} and Remark \ref{R2}, we conclude that the solution of {\rm (\ref{e8})} is $U$-recurrent. If $\alpha_3=0$, then the stochastic system {\rm (\ref{e8})} reduces to a geometric Brownian motion. It is known that, if $\alpha_1<\frac{\alpha_2}{2}$, then $\lim_{t\rightarrow\infty} X(t)=0$ a.s. {\rm (see {\rm\cite{Bernt:2005}})}. So the solution is recurrent with respect to the domain $U$. In this case, indeed, we can take $V_3(x)=\frac{x^2}{2}$ to obtain that
\begin{equation}\label{e9}
\mathcal LV_3(x)=\left(\alpha_1-\frac{\alpha_2}{2}\right) x^2,
\end{equation}
and the required assertion follows from Theorem \ref{T1}. If $\alpha_1>\frac{\alpha_2}{2}, \alpha_2>0$, we assume without loss of generality that $x_0>\delta$. Let $\sigma_n$ be the exit time sequence of $\{X(t)\}$ from $(\delta, 2^n \delta)$, where $n\in \nn$ and $x_0<2^n \delta. $  Using It\^o's formula to $\log(x)$ and Fatou's lemma, we have
$$
(\alpha_1-\frac{\alpha_2}{2})\mathbb {E}^{x_0}(\sigma_n)\leq(\alpha_1-\frac{\alpha_2}{2})\liminf_{t\rightarrow\infty} \mathbb {E}^{x_0}(\sigma_n\wedge t)\leq\log (2^n \delta)-\log(x_0).
$$
This means that $\mathbf E^{x_0}(\sigma_n)<\infty$. Furthermore, by It\^o's formula, we obtain
$$
\mathbb {E}^{x_0}(X(\sigma_n)^{r})=(x_0)^{r}, \q r=1-\frac{2\alpha_1}{\alpha_2}<0,
$$
which gives that
$$
p_n\delta^r+(1-p_n)(2^n\delta)^{r}=x_0^r, \q p_n=\mathbb {P}^{x_0}(X(\sigma_n)=\delta).
$$
Hence it follows that
$$
\lim_{n\rightarrow\infty}p_n=\mathbb {P}^{x_0}(\tau_U<\infty)=\left(\frac{x}{\delta}\right)^{r}<1.
$$
This also implies that  the solution of {\rm (\ref{e8})} is not $U$-recurrent when $\alpha_3=0$ and $\alpha_1>\alpha_2/2, \alpha_2>0$.
\end{example}
\par
In order to deal with multi-dimensional stochastic systems, one may adopt ideas and techniques in Example 3.1 and give the following conditions: Suppose that $U$ is a bounded open domain containing 0, $f$ and $\sigma$ satisfy the local Lipschitz condition (\ref{e2}). If there exists a constant $\alpha\in \rr$ and a function $\mu$ in class $\mathcal K$ such that
\begin{eqnarray*}
& &\frac{2x^T f(t,x)+\mbox{tr}(A(t,x))}{1+|x|^2}-\frac{2 x^T A(t,x) x}{(1+|x|^2)^2}\leq\alpha, \forall (t, x)\in \rr_+\times \rr^n, \\
&& \frac{2x^T f(t,x)+\mbox{tr}(A(t,x))}{|x|^2}-\frac{2 x^T A(t,x) x}{|x|^4}\leq -\mu(|x|), \forall (t, x)\in \rr_+\times U^c,
\end{eqnarray*}
where $A(t,x)=(a_{ij}(t,x))=\sigma(t,x)\sigma^T(t,x)$, then the solution of the equation (\ref{e1}) is recurrent relative to the domain $U$.
\begin{remark} A simple condition for recurrence relative to a bounded or unbounded domain $U$ is given by Khasminskii {\rm{(see \cite{Hasminskii:2012}, Theorem 3.9, p. 89)}}. More specifically, if the solution $X(t)$ of the equation {\rm (\ref{e1})} is regular and there exists in $\rr_+\times U^c$ a nonnegative function $V \in \mathbf C^{1,2}$, such that
\begin{equation} \label {e10}
\mathcal LV(s, x)\leq -\alpha(s),
\end{equation}
where $\alpha(s)\ge 0$ is a function for which
\begin{equation}
\beta(t)=\int_0^t \alpha(s)ds\rightarrow \infty, \; \mbox{as}\;\; t\rightarrow\infty,
\end{equation}
then $X(t)$ is $U$-recurrent.
\par
Note that the regularity of the solution $X(t)$ is needed in this criterion. Besides, if $U$ is a bounded domain, it is not easy to select an appropriate auxiliary  function $V\in \mathbf C^{1,2}$ in $\rr_+\times U^c$ with the property of {\rm (\ref{e10})} especially for those stochastic systems whose coefficients are independent of $t$.
\end{remark}
\par
If the domain $U^c$ is bounded, it will be more convenient to construct a nonnegative function $V$ satisfying (\ref{e10}) under certain assumptions. For example, if there exists
a pair of $a_{ii}(t,x)$ and $f_i(t, x)$ such that
\begin{equation} \label{e10-1}
0<a_0\leq a_{ii}(t, x), \q f_i(t,x)<b_0 \; \; \mbox{or}\;\; f_i(t,x)>b_0,\q b_0\in \rr,
\end{equation}
for any $(t,x)\in \rr_+\times U^c$, it has been shown in \cite {Hasminskii:2012} that a nonnegative function $V\in \mathbf C^{1,2}$ exists and satisfies
\begin{equation}
\mathcal LV(t,x)\leq -c e^{\gamma t}, \q c, \gamma>0.
\end{equation}
Hence the solution $X(t)$ of (\ref{e1}) is $U$-recurrent if the regularity holds.

Indeed, in this case, we also can construct a nonnegative function such that (\ref{e10}) holds under more general conditions. We shall need the following assumptions. For any $R>0$, there exists a pair of $a_{ii}(t,x)$ and $f_i(t, x)$ and positive constants $c_R, \hat c_R$ such that
\begin{equation}\label{e11}
a_{ii}(t,x) c_R-f_i(t,x)\ge \hat c_R, \q \mbox{if} \;\; |x|\leq R, t\in \rr_+,
\end{equation}
or
\begin{equation}\label{e12}
a_{ii}(t,x) c_R+f_i(t,x)\ge \hat c_R, \q \mbox{if} \;\; |x|\leq R, t\in \rr_+.
\end{equation}

Since $U^c$ is bounded, there exists a sufficiently large  $R>0$ so that $U^c\subset U_R=\{x: |x|<R\}$. Let $V(t, x)=K-(x_i-2R)^{2m}$, where the constants $K$ and $m\in \nn$ will be specified later. Assume that (\ref{e11}) holds. By It\^o's formula, for any $(t,x)\in \rr_+\times U_R$, it follows that
\begin{eqnarray*}
\mathcal LV(t,x)&=&-2m(x_i-2R)^{2m-1}f_i(t,x)-m(2m-1)(x_i-2R)^{2m-2}a_{ii}(t,x)\\
&= & -2m (x_i-2R)^{2m-2}\left(\frac{2m-1}{2}a_{ii}(t,x)+(x_i-2R)f_i(t,x)\right)\\
&\leq & -2m (x_i-2R)^{2m-2}\Big(3Rc_R a_{ii}(t,x)+(x_i-2R)f_i(t,x)\Big)\\
&\leq & -2m (x_i-2R)^{2m-2} R\hat c_R,
\end{eqnarray*}
if we choose an $m\in \nn $ by the condition $(2m-1)\ge 6R c_R$. We then take a sufficiently large positive constant $K\ge (3R)^{2m}$ so that $V(t,x)\ge 0$. It is easy to get that
$$
\mathcal LV(t,x)\leq -2m R \hat c_R(x_i-2R)^{2m-2}\leq -2m \hat c_R R^{2m-1}.
$$
If (\ref{e12}) holds, we can take $V(t,x)=\exp(\alpha R)-\exp(\alpha x_i)$, where $ \alpha$ is  a suitable constant. By It\^o's formula, we obtain that
\begin{eqnarray*}
\mathcal LV(t,x)&=&-\alpha \exp(\alpha x_i)f_i(t,x)-\frac{\alpha^2}{2} \exp(\alpha x_i)a_{ii}(t,x)\\
&= &- \alpha \exp(\alpha x_i)(\frac{\alpha}{2} a_{ii}(t,x)+f_i(t,x))\\
&\leq & -\sqrt{2c_R}\exp(\sqrt{2c_R}x_i)\hat c_R\\
&\leq&  -\sqrt{2c_R}\exp(-\sqrt{2c_R}R)\hat c_R,
\end{eqnarray*}
if we take $\alpha=\sqrt{2c_R}$. In terms of Theorem 3.9 in \cite{Hasminskii:2012}, we have proved the following theorem.
\begin{theorem} \label{T2}Under the conditions of {\rm{(\ref{e11})}} or {\rm{(\ref{e12})}}, if the domain $U^c$ is bounded, then the solution of {\rm (\ref{e1})} is recurrent relative to the domain $U$, provided that the regularity of {\rm(\ref{e1})} is satisfied.
\end{theorem}

In particular, if there exist two positive constants $r_1, r_2$ such that $U^c\subset \{r_1<|x|<r_2\}$,  Theorem 3.2 means that the first exit time of the solution $X(t)$ of (\ref{e1}) from the domain $U^c$ is finite with probability one, if $X(t)$ is regular and either (\ref{e11}) or (\ref{e12}) holds. In addition, it is easy to verify that either {\rm{(\ref{e11})}} or {\rm{(\ref{e12})}} will be satisfied, if there exists a pair of $a_{ii}(t,x)$ and $f_i(t, x)$ with the property of (\ref{e10-1}).

\par
\begin{theorem} \label{T3} Assume that {\rm (\ref{e2})} holds for the continuous coefficients $f$ and $\sigma$ in $\rr_+\times \rr^n$, and the domain $U$ is the same as that in Theorem 3.1. Suppose moreover that there exists a nonnegative function $V\in \mathbf C^{2}$ in $\rr^n\setminus \{0\}$ such that for any $a>0$,
\begin{eqnarray}\label{e12-1}
V(x)=\Phi(|x|), \q  \mathcal LV(x)\leq 0 \;\; \mbox{as}\;\; |x|\ge a
\end{eqnarray}
and
\begin{eqnarray}\label{e13}
\Phi'(r)<0, \;\; \lim_{r\rightarrow \infty } \Phi(r)=0\;\; \mbox{as}\;\; r\ge a.
\end{eqnarray}
If there exists a pair of $a_{ii}(t,x)$ and $b_i(t,x)$ such that either {\rm{(\ref{e11})}} or {\rm{(\ref{e12})}} holds for any $R>0$, then
the solution $X(t)$ of {\rm(\ref{e1})} is not recurrent relative to the domain $U$ for any initial value $x_0\not\in [U]$.
\end{theorem}
\par
{\bf Proof.} Since $f$ and $\sigma$ satisfy the local Lipschitz condition (\ref{e2}), there exists a unique maximum local solution $X(t)$ to the equation (\ref{e1}). Let $\tau_n$ be the first exist time of the process $X(t)$ from the domain $\{x\in \rr^n: |x|<n\}$. It is clear that the sequence of stopping times
$\{\tau_n\}$ is monotonically increasing, and its limit, denoted by $\tau$, is called the explosion time of $X(t)$. If $\mathbb{P}^{x_0}(\tau<\infty)>0$, then the solution $X(t)$ is not regular, hence it is not recurrent relative to the domain
$U$ either. It remains to prove the assertion of the theorem when $\mathbb {P}^{x_0}(\tau=\infty)=1$.  Let $\tau_U$ be the first time that the path of $x(t)$ reaches the boundary of $U$. Since $U$ is bounded, we can find an $\alpha>0$ satisfying $|x_0|>\alpha$ so that $U\subset\{x\in\rr^n: |x|<\alpha\}$. We now choose a sufficiently large $\beta>\alpha>0$ such that $x_0\in D_{\alpha, \beta}=\{x\in \rr^n:
\alpha<|x|<\beta\}$. Let $\tau_{\alpha,\beta}$ denote the exit time from the shell $D_{\alpha, \beta}$. By Theorem \ref{T2}, it follows that
$\tau_{\alpha,\beta}<\infty$ a.s. By It\^o's formula, we have
$$
\mathbb{ E}^{x_0}\Phi(|X(\tau_{\alpha, \beta}\wedge t)|)=\Phi(|x_0|)+\mathbb {E}^{x_0}\int_0^{\tau_{\alpha, \beta}\wedge t}\mathcal LV(X(s))ds\leq \Phi(|x_0|).
$$
Taking $t\rightarrow\infty$ and using Fatou's lemma, we get
\begin{equation*}
\mathbb{ E}^{x_0}\Phi(|X(\tau_{\alpha, \beta})|)\leq \Phi(|x_0|),
\end{equation*}
which implies that
\begin{eqnarray*}
\mathbb {P}^{x_0}(|X(\tau_{\alpha, \beta})|=\alpha)\Phi(\alpha)+\mathbb {P}^{x_0}(|X(\tau_{\alpha, \beta})|=\beta)\Phi(\beta)\leq \Phi(|x_0|).
\end{eqnarray*}
Hence, taking $\beta\rightarrow\infty $ and using (\ref{e13}), we have
$$
\lim_{\beta\rightarrow\infty}\mathbb {P}^{x_0}(|X(\tau_{\alpha, \beta})|=\alpha)\leq \frac{\Phi(|x_0|)}{\Phi(\alpha)}.
$$
By this and (\ref{e13}) again, one easily see that
$$
\mathbb {P}^{x_0}(\tau_U<\infty)\leq \lim_{\beta\rightarrow\infty}\mathbb {P}^{x_0}(|X(\tau_{\alpha, \beta})|=\alpha)\leq \frac{\Phi(|x_0|)}{\Phi(\alpha)}<1.
$$
This gives the assertion. \hfill $\Box$
\begin{example} One-dimensional O-U process is described by
\begin{equation}
d X_t=\mu X_t dt +\sigma d B_t, \q \sigma> 0, \quad X_0=x_0.
\end{equation}
This equation has a unique solution and the solution is given by
$$
X_t=e^{\mu t}x_0+\sigma \int_0^t e^{\mu(t-s)}dB_s.
$$
Let $\delta>0$ and $U=(-\delta, \delta)$. If $\mu\leq 0$, we take
$$
V(x)=\int_0^{|x|} \exp\big[-\frac{\mu}{\sigma^2} y^2\big] dy.
$$
It is easy to derive that
$\mathcal LV(x)=0$. Noting that $|\frac{\partial V}{\partial x}\sigma|^2$ is a $\mathcal K_\infty$ class function. Hence, by Theorem \ref{T1} the solution $X_t$ is recurrent relative to $U$ for any initial value $x_0\not\in [-\delta,\delta]$. If $\mu>0$, we take
$$
V(x)=\int_{|x|}^\infty \exp\big[-\frac{\mu}{\sigma^2} y^2\big] dy.
$$
It is clear that the conditions {\rm (\ref{e11}), (\ref{e12-1})} and {\rm (\ref{e13})} hold, and the solution $X_t$ is not recurrent with respect to the domain
$U$ for any $x_0\not\in [-\delta,\delta]$ according to Theorem {\rm\ref{T3}}. This example, indeed, shows that $\mu\leq 0$ is the sufficient and necessary condition for the recurrence of $U$.
\end{example}

Theorem \ref{T3} requires that there exists a function $\Phi$ such that both (\ref{e12-1}) and (\ref{e13}) are satisfied. To this end, inspired by the assumptions of Theorem 9.1.1 in \cite{Friedman:1975}, we shall construct a suitable function $\Phi$ under some conditions.

We impose the following assumptions:
\par
(A$_1$) The matrix $(a_{ij}(t,x))$ is positive definite for any $(t,x)\in \rr_+\times \rr^n$.
\par
(A$_2$) For any $a>0$,
$$
S(t,x)\ge \theta(|x|), \q \mbox{if} \;\; |x|\ge a,
$$
where $\theta: \rr_+\rightarrow \rr$ is a continuous function such that
\begin{equation}\label{e14}
\int_a^\infty \exp\big(-\int_a^r\frac{\theta(s)}{s}ds\big)dr<\infty;
\end{equation}
and
$$
S(t,x)=\frac{B(t, x)+C(t,x)-A(t,x)}{A(t, x)},
$$
\begin{eqnarray*}
A(t,x)&=& \frac{1}{|x|^2}\sum_{i,j=1}^n a_{ij}(t,x)x_ix_j, \q B(t,x)=\sum_{i=1}^n a_{ii}(t,x),\\
C(t,x)&=& 2\sum_{i=1}^n x_i f_i(t,x).\\
\end{eqnarray*}
Notice that (\ref{e14}) is fulfilled  for any of the functions
$$
\theta(s)=k_1 s^2, k_1>0, \q \theta(s)=k_2 s, k_2>0, \q \theta(s)=d>1.
$$
An application of It\^o's formula to $V(x)=\Phi(|x|)$ yields that
\begin{eqnarray*}
\mathcal LV(x)&=& \frac{1}{2}A(t,x)\left(\Phi''(|x|)+\frac{\Phi'(|x|)}{|x|} S(t,x)\right)\\
&\leq & \frac{1}{2} A(t,x)\left(\Phi''(|x|)+\frac{\Phi'(|x|)}{|x|} \theta(|x|)\right)
\end{eqnarray*}
if $\Phi'(|x|)<0$ for $|x|\ge a$. From (\ref{e14}), we can take
$$
\Phi(r)=\int_r^\infty \exp (-\int_a^t\frac{\theta(s)}{s}ds)dt, \q r\ge a.
$$
It is easy to verify that the conditions (\ref {e12-1}) and (\ref{e13}) are satisfied, since
$$
\Phi''(r)+\frac{\Phi'(r)}{r} \theta(r)=0, \q r\ge a.
$$
Hence, under the assumptions of (A$_1$) and (A$_2$), if either {\rm{(\ref{e11})}} or {\rm{(\ref{e12})}} holds, then the solution of {\rm(\ref{e1})} is not recurrent relative to the domain $U$ as in Theorem \ref{T3} for any initial value $x_0\not\in [U]$

If the condition (\ref{e6-3}) in Theorem \ref{T1} is replaced by a somewhat more stringent condition, the following theorem is useful for estimating the expectation of the stochastic residence time.

\begin{theorem} \label {T4}Given a domain $U$ as in Theorem \ref{T1}. Assume that {\rm (\ref{e2})} holds for the continuous coefficients $f$ and $\sigma$ in $\rr_+\times \rr^n$. If there exists a nonnegative function
$V\in \mathbf C^{1,2}$ such that
\begin{eqnarray}
& &V(t,x)\ge \mu_1(|x|), \q \mu_1\in \mathcal K, \q \label{e15}\\
& &\mathcal LV(t,x)\leq -\lambda V(t, x), \lambda>0, \label{e16}
\end{eqnarray}
then the solution $X(t)$ of the equation {\rm (\ref{e1})} is $U$-recurrent and satisfies
\begin{eqnarray}
\mathbb {E}^{x_0} \left( e^{\lambda \tau_U}\right)\leq \frac{V(0, x_0)}{\inf_{t\ge 0, x\in \partial U}V(t, x)},\label {e19},
\end{eqnarray}
where $T>0$ and $x_0\in U^c$.
\end{theorem}
\par
{\bf Proof.} By (\ref{e2}), (\ref{e15}) and (\ref{e16}), it follows from Theorem \ref{l1} that the regularity holds for the equation (\ref{e1}) with arbitrary given initial value $X(0)=x_0\in \rr^n$. Moreover, employing (\ref{e15}), (\ref{e16}) again and Theorem \ref{T1}, we conclude that the solution of the equation (\ref{e1}) is $U$-recurrent, that is,
$\mathbb {P}^{x_0}(\tau_U<\infty)=1$ for any $x_0\in U^c$. Let $\tau_n$ be the first exist time of the process $X(t)$ from the domain $\{x\in \rr^n: |x|<n\}$. Applying It\^o's formula, for any $T>0$ we have
\begin{equation}
\mathbb {E}^{x_0}\left( e^{\lambda( \tau_U\wedge\tau_n\wedge T)}V(\tau_U\wedge\tau_n\wedge T, X(\tau_U\wedge\tau_n\wedge T))\right)\leq V(0, x_0).
\end{equation}
Since the solution of the equation (\ref{e1}) is regular, Fatou's lemma yields
\begin{equation}\label{e17}
\mathbb {E}^{x_0}\left( e^{\lambda \tau_U}V(\tau_U, X(\tau_U))\right)\leq V(0, x_0).
\end{equation}
Note that
$$
\inf_{t\ge 0, x\in \partial U}V(t, x)\ge \mu_1(r)>0, \q r=\min_{x\in \partial U}|x|>0,
$$
due to (\ref{e15}). It then follows from (\ref{e17}) that
$$
\inf_{t\ge 0, x\in \partial U}V(t, x) \mathbb {E}^{x_0}\left (e^{\lambda \tau_U}\right)\leq V(0, x_0),
$$
which gives the required assertions.  \hfill $\Box$

\par
By (\ref{e19}) of Theorem \ref{T4}, it is easy to get
\begin{equation}\label{e20}
\mathbb {P}^{x_0}(\tau_U>T)\leq \frac{V(0, x_0)}{e^{\lambda T}\inf_{t\ge 0, x\in \partial U}V(t, x)},
\end{equation}
via Chebyshev's inequality. This inequality can be used to deal with the domain control problems, just as stated in Section 1. In particular, if the bounded domain $U$ is an open neighborhood $U_\delta$, and the inequalities (\ref{e15}) and (\ref{e16}) are replaced by the following equalities
\begin{eqnarray}
V(t,x)&=&\mu_1(|x|), \q \mu_1\in \mathcal K, \\
\mathcal LV(t,x)&=& -\lambda V(t, x), \q \lambda>0,
\end{eqnarray}
then it is obvious that
$$
\mathbb {E}^{x_0} \left( e^{\lambda \tau_U}\right)=\frac{V(0, x_0)}{\mu_1(\delta )}=\frac{\mu_1(|x_0|)}{\mu_1(\delta )}.
$$
\par

Let $D$ be a domain in $\rr^n$. We shall say that the partial differential operator $\mathcal L$ is {\sl uniformly elliptic} in $\rr_+\times D$, if there exists a positive constant $\mu$ such that
 $$
 \sum_{i,j=1}^n a_{ij}(t,x)\xi_i\xi_j\ge \mu |\xi|^2,
 $$
  for all $t\in \rr_+, x\in D, \xi\in \rr^n.$ That is, all the eigenvalues of the symmetric matrix $(a_{ij}(t,x))_{i,j=1}^n$ have a positive lower bound.

\par
The following theorem will give a connection between $\mathbb {E}^{x}(\tau_U)$ and a Dirichlet problem under suitable assumptions. The proof uses an approach similar to that used for proving Theorem 3.11 in \cite{Hasminskii:2012}.

\begin{theorem}\label{T5}
Let $U$ be an open ball $B_\delta=\{x\in \rr^n: |x|<\delta\}, 0< \delta<\infty$. Assume that the coefficients $f$ and $\sigma$ are independent of $t$ and satisfy the local Lipschitz condition {\rm (\ref{e2})}. Also, assume that there exists a nonnegative function $V\in \mathbf C^{1, 2}$ such that {\rm(\ref{e6-1})} and {\rm(\ref{e6-3})} hold in $\rr_+\times B_\delta^c$. Then $\mathbb {E}^{x_0}(\tau_U)$ exists for any initial value $x_0\not\in U$ of the equation {\rm(\ref{e1})} and is a solution of the following Dirichlet problem
\begin{eqnarray}
&&\mathcal Lu(x)=-1, \q x\in [U]^c, \\
&&u(x)|\partial U=0,
\end{eqnarray}
provided that the partial differential operator $\mathcal L$ is uniformly elliptic in $U^c$.

\end{theorem}

\par
{\bf Proof.} Under the conditions of the theorem, it follows from Theorem \ref{l1} and Remark \ref{r1} that  the regularity of the equation (\ref{e1}) holds. From Theorem \ref{T1} and Remark \ref{R2}, we can conclude that
$\mathbb {E}^{x_0}(\tau_U)<\infty$ for any given initial value $x_0\not\in U$. Let $\tau_U^n$ be the first exit time from the open domain $[U]^c\cap \{|x|<n\}$. Due to the regularity of the solution, it is obvious that
\begin{equation}\label{e23}
\mathbb {E}^{x_0}(\tau_U^n)\rightarrow \mathbb {E}^{x_0}(\tau_U), \q n\rightarrow \infty,
\end{equation}
by the monotone convergence theorem. Note that the coefficients $f$ and $\sigma$ satisfy the Lipschitz condion in the domain $U^c\cap \{|x|\leq n\}$ and the diffusion matrix is uniformly elliptic, it then follows from the theory of partial differential equation (see e.g. Friedman \cite{Friedman:1975}, Theorem 2.4, p.134) that the following Dirichlet problem
\begin{eqnarray}
&&\mathcal Lu_n(x)=-1,\q x\in [U]^c\cap \{|x|<n\}, \label{e21}\\
&& u_n(x)\big|\partial U=0, \q u_n(x)\big|\{|x|=n\}=0,\label{e22}
\end{eqnarray}
admits a unique solution. By means of It\^o's formula, one easily sees that $\mathbb {E}^{x}(\tau_U^n)$ is the unique solution of (\ref{e21}) and (\ref{e22}) by virtue of the definition of $\tau_U^n$. We now take a sufficiently large $n_0\in \mathbb{N}$ such that $[U]\subset \{|x|\leq n_0\}$. For any $n\ge n_0$, let

$$
\hat u_n(x)=\left\{\begin{array}{ll} \mathbb {E}^{x}(\tau_U^n), &\q\mbox{if}\;\; x\in U^c\bigcap \{|x|<n\},\\
0, &\q \mbox{if}\;\; x\in \{|x|\ge n\}.
\end{array}
\right.
$$
We now define
$$
v_n(x)=\hat u_{n+1}(x)-\hat u_n(x), x\in [U]^c, n\ge n_0.
$$
It is obvious that $v_n(x)\ge 0$ and satisfies
$$
\mathcal Lv_n(x)\leq 0, \q \forall x\in [U]^c, n\ge n_0.
$$
Note that for any given initial value $x_0\in [U]^c$, there exists a sufficiently large $n_1\ge n_0$ such that $x_0\in [U]^c\cap \{|x|<n\}$ for each $n\ge n_1$. In this case, it is easy to know that
\begin{eqnarray}
& &\mathcal Lv_n(x_0)=0, \q  n\ge n_1, \\
& & \mathbb {E}^{x_0}(\tau_U):=u(x_0)=\hat u_{n_1}(x_0)+\sum_{n=n_1}^\infty v_n(x_0),
\end{eqnarray}
due to (\ref{e23}). By Harnack's second theorem of harmonic functions, it can be shown that $\sum_{n=n_1}^\infty v_n(\cdot)$ is harmonic at $x_0$ (see, for example, \cite{Gilbarg:2001}, p. 21-22), if necessary, we can take a small neighborhood of $x_0$, denoted by $U_\delta(x_0)$, and consider $v_n(x)$ for all $x\in U_{\delta}(x_0)$. This, together with (\ref{e21}), implies that
$\mathcal Lu(x_0)=-1$. Since $x_0$ is arbitrary, then $\mathcal Lu(x)=-1$ holds for any $x\in [U]^c$. By the definition of $\tau_U$ and the property of $\mathbb{E}^{x}(\tau_U)$, we thus have $u(x)\big|_{\partial U}=0$ as required.  \hfill $\Box$

\begin{example} Consider a scalar stochastic system in the form
\begin{equation}\label{e23-1}
dX(t)=f(X(t))dt+dB_t, \q X(0)=x_0\in (1, \infty)\cup(-\infty, -1).
\end{equation}
Assume that the drift term $f(x)=-x^m$ for some odd number $m\in \mathbb{N}$. In this case, by taking $V(x)=x^2$, it is easy to derive that
$$
\mathcal LV(x)=-2x^{m+1}+1.
$$
If $U=(-1, 1)$, it then from Theorem {\rm\ref{l1}}, Theorem {\rm{\ref{T1}}} and Remark {\rm\ref{R2}} that the unique solution of {\rm (\ref{e23-1})} is recurrent relative to the interval $U$ for any given initial value
$x_0\not\in(-1, 1)$. To be precise, we have
\begin{equation} \label{e23-2}
\mathbb {E}^{x_0}(\tau_U)\leq \frac{x_0^2-1}{2-1}=x_0^2-1.
\end{equation}
For the system {\rm (\ref{e23-1})}, we can conclude from Theorem {\rm \ref{T5}} that $\mathbb {E}^{x}(\tau_U):=\tau(x)$ satisfies
\begin{equation}\label{e24}
\left\{\begin{array}{ll}
-x^m\frac{d\tau(x)}{dx}+\frac{1}{2}\frac{d^2\tau(x)}{dx^2}=-1, &\q x\not\in [-1, 1],\\
 \tau(1)=\tau(-1)=0. &
\end{array}
\right.
\end{equation}
Even in such a simpler situation, it is very difficult to find an exact solution for $\tau(x)$. Therefore, we shall solve {\rm(\ref{e24})} numerically with the help of the R software for various initial values of $X(0)$. Table 1. shows the computational results for the Dirichlet problem of {\rm(\ref{e24})} relative to the stochastic system {\rm(\ref{e23-1})} with $m=1$ and $m=3$  by setting an upper bound of $\tau(1+h)$ according to {\rm (\ref{e6-4})} of Remark {\rm \ref{R2}}, where $h=(3-1)/10000$ is the step size. We need to search the best value of $\tau(1+h)$ via {\rm(\ref{e24})} such that the following quantity
$$
\frac{(2+2h)\tau(1+h)-\tau(1+2h)-2h^2}{1+2h}
$$
achieves the least value since it is an approximation of $\tau(1)=0$.
\begin{table}[!htbp]
\caption{Numerical results on the mean of $\tau_U$ for various initial values}
\centering
\begin{tabular}{lccl}
  \hline
 $f(x)$ &\q\q\q \q$x_0$\q\q\q \q \q &$\tau(x_0)=\mathbb {E}^{x_0}(\tau_U)$  & \q $\hat\tau(x_0)$ (by (\ref{e23-2}))  \\
  \hline
  & 1.5 &\q  0.4812779 & \q\q 1.25\\
  $-x$ & 2.0 &\q  0.7956221 & \q\q 3.0 \\
  & 2.5 & \q 1.0306462 & \q \q5.25\\
  & 3.0 & \q  1.2193076 & \q \q8.0\\
  \hline
  & 1.5 &  \q 0.3231637 & \q\q 1.25\\
  $-x^3$ & 2.0 &\q  0.4235235 & \q \q 3.0\\
  &2.5 & \q   0.4690177 & \q \q 5.25\\
  & 3.0 &\q 0.4935798 & \q \q 8.0\\
  \hline
  \end{tabular}
  \end{table}
\end{example}
\section {Domain Aiming Control}

In this section, we first consider a controlled stochastic system of the form
\begin{equation}\label{e25}
dX(t)=(g(t,X(t))+u)dt+\sigma(t, X(t))dB_t, \q X(0)=x_0\in \rr^n,
\end{equation}
where $\sigma:\rr_+\times \rr^n\rightarrow \rr^m$ satisfies the local Lipschitz condition. Assume that there exists a local Lipschitz continuous function $f: \rr^+\times \rr^n \rightarrow \rr^n$ and a Lyapunov function $V\in\mathbf C^{1,2}$ such that either conditions (\ref{e6-1}) and (\ref{e6-3}) of Theorem {\ref{T1} hold or conditions (\ref{e15}) and (\ref{e16}) of Theorem \ref{T4} are satisfied in $\rr_+\times U^c$ for an open bounded domain $U$ containing 0 in its interior. Without loss of generality, we assume that $U$ is an open ball $B_\delta$ for simplicity, where $\delta$ is a positive constant.
\par
Given a pair $(T, p)$ with $T>0$ and $0<p<1$, we are interested in choosing a Markovian control law $u(t)$, by which the state $X(t)$ of the resulting closed-loop stochastic system, originating from $x_0$ outside the target domain $B_\delta$, will reach the boundary of $B_\delta$ during period $T$ with probability $p$.
\begin{definition} Given a triple $(T, p, U)$, if there exists an admissible control $u(t)=u(t, X(t))$ such that the resulting stochastic system of the controlled stochastic system {\rm (\ref{e0})} has
the property of
$$
\mathbb{P}^{x_0}\left(\exists t\in [0, T], X(t,x_0)\in \partial U\right)\ge p,
$$
for any given initial value $X(0)=x_0\in [U]^c$, then the system {\rm (\ref{e0})} is said to be residence probability controllable in the domain $[U]^c$. In this case, $u(t)$ is called a domain aiming controller of $U$.
\end{definition}

\begin{theorem} \label {T6} For any given pair $(T, p)$, under the above assumptions, the controlled stochastic system {\rm(\ref{e25})} is residence probability controllable with respect to the domain $[B_\delta]^c$, if there exists a nonnegative function $V\in \mathbf C^{1,2}$ in $\rr_+\times B_{\delta}^c$ such that the initial value $x_0\in B_{\delta}^c$ either satisfies
\begin{eqnarray}
V(0, x_0)\leq \inf_{t\ge 0, |x|=\delta}V(t,x)+T(1-p)\mu(\delta),
\end{eqnarray}
under the conditions of {\rm (\ref{e6-1})} and {\rm (\ref{e6-3})} in Theorem {\rm\ref{T1}}, or satisfies
\begin{equation}
V(0, x_0)\leq e^{\lambda T}(1-p)\mu_1(\delta),
\end{equation}
under the conditions of {\rm(\ref{e15})} and {\rm(\ref{e16})} in Theorem {\rm \ref{T4}}.
\end{theorem}
\par
{\bf Proof.} Under the assumptions of the theorem, for each given pair $(T, p)$, we can choose a Markovian control law given by
 $$
 u=-g(t, X(t))+f(t, X(t)).
 $$
 According to Theorem \ref{T1}, Remark \ref{R2} and Theorem \ref{T4}, we can conclude that the corresponding closed-loop stochastic system is regular and its solution is recurrent with respect to the domain $B_{\delta}$. Let $\tau$ be the residence time in the domain $[B_\delta]^c$ for any given initial value $x_0$ therein. We thus have
 $$
 \mathbb {E}^{x_0}(\tau)\leq \frac{V(0, x_0)+\int_0^\infty \nu(t)dt-\inf_{t\ge 0, |x|=\delta}V(t, x)}{\mu(\delta)},
  $$
 if the conditions (\ref{e6-1}) and (\ref{e6-3}) are fulfilled. By using Chebyshev's inequality, we further have
 $$
 \mathbb {P}^{x_0}(\tau> T)\leq \frac{\mathbb {E}^{x_0}(\tau)}{T}\leq \frac{ V(0, x_0)+\int_0^\infty \nu(t)dt-\inf_{t\ge 0, |x|=\delta}V(t, x)}{T\cdot \mu(\delta)},
 $$
 which implies that $ \mathbb {P}^{x_0}(\tau\leq T)\ge p $ holds, if the initial value $x_0$ satisfies
 $$
 \frac{V(0, x_0)+\int_0^\infty \nu(t)dt-\inf_{t\ge 0, |x|=\delta}V(t, x)}{\mu(\delta)}\leq T(1-p).
 $$
 Under the conditions (\ref{e15}) and (\ref{e16}) of Theorem \ref{T4}, it follows from (\ref{e20}) that $\mathbb {P}^{x_0}(\tau\leq T)\ge p $ still holds, provided that
 $$
\frac{V(0, x_0)}{\mu_1(\delta)}\leq e^{\lambda T}(1-p)
 $$
 is satisfied. The proof is complete. \hfill $\Box$
\begin{example} Consider the following controlled stochastic system
\begin{equation} \label{e26}
dX(t)=(g(t, X(t))+u)dt+\hat\sigma(t, X(t))\cdot X(t) dB_t, \q X(0)=x_0\in \rr^n,
\end{equation}
where $\hat\sigma: \rr^{+}\times \rr^n\rightarrow \rr$ is local Lipschitz continuous and satisfies
\begin{equation}\label{e27}
|\hat\sigma(t,x)|^2\ge \alpha(1+|x|^2), \q \alpha>0.
\end{equation}
Let
$$
u=u(t,x)=-g(t,x)-\hat\sigma(t,x)^2 x.
$$
Now, choosing a Lyapunov function $V$ defined on $\rr^n$ by $V(x)=\frac{1}{2}|x|^2$, one can prove that the closed-loop system
$$
dX(t)=-\hat\sigma(t, X(t))^2\cdot X(t)dt+\hat\sigma(t, X(t))\cdot X(t) dB_t, \q X(0)=x_0\in \rr^n,
$$
is regular and recurrent with respect to the open ball $B_1=\{|x|<1\}$ for any given initial value $x_0\not\in B_1$ from 
Theorem {\rm \ref{l1}} and Theorem {\rm\ref{T1}}, since
$$
\mathcal LV(x)\leq -\frac{1}{2}\alpha |x|^2(1+|x|^2).
$$
By Theorem {\rm\ref{T6}}, for any given pair $(T, p)$, if the initial value $x_0$ and $\alpha$ satisfy
$$
|x_0|^2\leq 1+2\alpha\cdot T\cdot (1-p),
$$
then the controlled stochastic system {\rm(\ref{e26})} is residence probability controllable outside the open ball $B_1$. Note that the condition {\rm(\ref{e27})} can be removed for the problem of residence probability control indeed. To see this, we can choose the same Lyapunov function $V$ as above and a control law
$$
u=u(t, x)=-g(t,x)-\frac{1}{2}\hat\sigma(t,x)^2 x-\frac{1}{2}x.
$$
Similarly, one can derive that
$$
\mathcal LV(x)=-\frac{1}{2}|x|^2=-V(x).
$$
By Theorem {\rm\ref{T6}} again, it is readily seen that if
$$
|x_0|^2\leq e^T (1-p),
$$
then the system {\rm(\ref{e26})} is also residence probability controllable outside $B_1$.
\end{example}

\par
In the case of linear controlled stochastic systems, Theorem \ref{T6} can be used to design an appropriate admissible control such that the corresponding closed-loop stochastic system is residence probability controllable under unrestrictive conditions.
\par
Let us consider a linear controlled stochastic system in the form
\begin{equation} \label{e28}
dX(t)=(AX(t)+Bu)dt+CdB_t, X(0)=x_0\in\rr^n,
\end{equation}
where $A$ is a square matrix of order $n$, while $B$ and $C$ are $n\times l$ and $n\times m$ matrices, respectively. If $n=l$ and $B$ is of full rank, we have:
\begin{theorem} \label{T7}
Assume that there exists a Hurwitz matrix $D$ such that the pair $(D, C)$ is completely disturbable, namely, $rank (C|DC|\cdots|D^{n-1}C)=n$. For any given pair $(T, p)$, the stochastic system {\rm(\ref{e28})} is residence probability controllable outside an open ball $B_\delta$ for any initial value $x_0 \in B_\delta^c$.
\end{theorem}
{\bf Proof.} Since $D$ is Hurwitz and $(D, C)$ is disturbable, then there exists a unique positive definite solution $M$ to the matrix equation
\begin{equation*}
MD+D^TM+ MCC^TM=0,
\end{equation*}
which is equivalent to
\begin{equation} \label{e29}
DM^{-1}+M^{-1}D^T+CC^T=0.
\end{equation}
Let $V(x)=\frac{1}{2}x^T M^{-1} x$ and define
$$
u=u(t,x)=B^{-1}\big(-A+D^T-\gamma I_n\big)x,
$$
where $I_n$ is the identity matrix and $\gamma>0$ is a sufficiently large constant that will be determined later. It is easy to deduce that
\begin{eqnarray*}
\mathcal LV(x)&=&\frac{1}{2}x^T(M^{-1}D^T+D M^{-1})x-\gamma x^TM^{-1}x+\frac{1}{2} \;\mbox{trace}\;(M^{-1}CC^T)\\
&=& -\frac{1}{2}x^T(CC^T)x-\gamma x^TM^{-1}x+\frac{1}{2} \;\mbox{trace}\;(M^{-1}CC^T)\\
&\leq& -\left(\frac{\lambda_1}{2}+\frac{\gamma}{\lambda_{\max}(M)}\right)|x|^2+\frac{1}{2\lambda_{\min}(M)}\sum_{i=1}^n \lambda_i, \label{e30}
\end{eqnarray*}
where $\lambda_i, i=1, 2,\cdots n$ are eigenvalues of $CC^T$, and $\lambda_1\ge 0$ is the smallest one. If we choose an appropriate $\gamma>0$ such that
$$
\left(\frac{\lambda_1}{2}+\frac{\gamma}{\lambda_{\max}(M)}\right)|\delta|^2-\frac{1}{2\lambda_{\min}(M)}\sum_{i=1}^n \lambda_i\triangleq a_{\gamma}>0,
$$
then from Theorem \ref{T1} and Remark \ref{R2}, it is known that the resulting closed-loop system is recurrent with respect to $B_\delta$. By Theorem \ref{T6}, furthermore, if the initial value $x_0$ satisfies
\begin{equation}
\frac{1}{2}x_0^T M^{-1}x_0\leq \frac{1}{2\lambda_{\max}(M)}\delta^2+T(1-p) a_\lambda,
\end{equation}
then the required conclusion follows. Indeed, for any given $T, p, \delta$ and $x_0$, we always can find a sufficiently large positive constant $\gamma$ satisfying (4.8), since $\lambda_{\min}(M)>0$ and $\lim_{\gamma\rightarrow\infty}a_\gamma\rightarrow\infty$. \hfill $\Box$
\begin{remark} If $B$ has full column rank with $l<n$, then by an suitable linear transform $ \bar X=P  X$, the system {\rm(\ref{e28})} can be reduced to the form
\begin{equation}\label{e31}
\left[ \begin{array}{c}
d\bar X_1(t) \\
d\bar X_2(t) \\
\end{array}
\right ]=\left[\begin{array}{cc}
\bar A_{11} & \bar A_{12}\\
\bar A_{21} & \bar A_{22}\\
\end{array}\right]\left[ \begin{array}{c}
\bar X_1(t)\\
\bar X_2(t)\\
\end{array}\right]dt +\left[\begin{array}{c}
\bar B_1\\
0\\
\end{array} \right] udt +\left[\begin{array}{c}
\bar C_1\\
\bar C_2\\
\end{array}\right]dB_t,
\end{equation}
where $P$ and $\bar B_1$ are $n\times n$ and $l\times l$ matrices with full rank. Furthermore, if $\bar A_{21}=0$ and $\bar C_2=0$, and the deterministic subsystem,
$d\bar X_2(t)=\bar A_{22}\bar X_2(t)dt$, is asymptotically stable in the sense that  $\lim_{t\rightarrow\infty}\bar X_{2}(t)=0$ for any initial value $\bar x_2(0)$, we can consider an alternative domain aiming control problem for the controlled stochastic system
$$
d\bar X_1(t)=(\bar A_{11} \bar X_1(t)+\bar A_{12} \bar X_2(t))dt+\bar B_1 udt+\bar C_1dB_t,
$$
for which Theorem {\rm\ref{T6}} can apply. It is obvious that, if the pair $(A, B)$ is nonreachable (see, for example {\rm\cite {Zak:2003}}, Theorem {\rm 3.4}, p.103), then there exists a similarity transform $\bar X=P X$ such that $\bar A_{21}=0$ in {\rm (\ref{e31})}. In addition, if $\mbox{Span}(C)\subset\mbox{Span}(B)$, then $\bar C_2=0$, where $\mbox{Span}(C)$ denotes the subspace of $\rr^n$ generated by all column vectors of $C$.
\end{remark}

\section{Conclusion}
In this paper, a problem of domain aiming control has been formulated and considered for controlled stochastic systems described by It\^o stochastic differential equations. The goal of control is to find a control law such that the trajectory of the close-loop system, originating from a initial value outside an open and bounded domain, can reach the domain during a specified duration with a certain probability. The approach adopted in this paper was based on probabilistic analysis of the stochastic residence time outside the domain. It is often difficult or even impossible to get the probability distribution of the residence time. Alternatively, Lyapunov-type conditions have been given for obtaining an upper bound to the expected residence time. If the target domain is nonrecurrent, it is generally infeasible to solve this issue by designing a desired controller. Accordingly, various criteria for the recurrence and non-recurrence relative to a bounded open domain or an unbounded domain have been provided by means of Lyapunov functions.

\par

Future research effort can be devoted to the problem of domain aiming control for different types of controlled stochastic nonlinear systems. For example,
the drift of the stochastic system is affine in the control input. It is also of interest to investigate the case that the diffusion of the stochastic system has a control input. Another direction that is of extremely important in theory is to give sufficient and necessary conditions for domain aiming controllability of stochastic systems.


\begin{thebibliography}{00}

\bibitem{Battilotti:2001}
S.~Battilotti, A unifying framework for the semiglobal stabilization of nonlinear uncertain systems via measurement feedback,
{\it IEEE Trans. Automat. Contr.}, 46(1), 3-16, 2001.


\bibitem {Battilotti:2003}

S.~Battilotti and A.~Santis, Stabilization in probability of nonlinear stochastic sytems with guaranteed region of attraction and target set,
{\it IEEE Trans. Automat. Contr.}, 48(9), 1585-1599, 2003.

\bibitem{Deng:2001}

H.~Deng, M.~Krsti\'{c}, and R.~J.~Williams, Stabilization of stochastic
nonlinear systems driven by noise of unknown covariance, {\it IEEE
Trans. Automat. Contr.}, 46(8): 1237¨C1253, 2001.


\bibitem {Friedman:1975}
A.~Friedman, Stochastic Differential Equations and Applications, Volume 1, Academic Press, Now York, 1975.

\bibitem{Garnell:1977}
P.~Garnell, Guided Weapon Control Systems, Pergamon Press, New York, 1977.

\bibitem{Gilbarg:2001}

D.~Gilbarg and N.~S.~Trudinger, Elliptic Partial Differential Equations of Second Order,
Springer-Verlag, Berlin, 2001.

\bibitem {Karatzas:1999}
I.~Karatzas and S. Shreve, Brownian Motion and Stochastic Calculus, 2nd Edition, Springer-Verlage, New York, 1999.


\bibitem{Khalil:2002}
H.~K.~Khalil, Nonlinear Systems, Third Edition, Prentice Hall, New Jersy,  2002.

\bibitem {Hasminskii:2012}
R.~Khasminskii,  Stochastic Stability of Differential
Equations, Second Edition, Springer, 2012.



\bibitem{Khoo2013}
S.~Khoo, J.~Yin, Z.~Man and X.~Yu,  Finite-time stabilization of stochastic nonlinear systems
in strict-feedback form, {\it Automatica,} 49(5): 1403-1410, 2013. 

\bibitem{Kim:1992}
S.~Kim, S.~M.~Meerkov and T.~Runolfsson, Aiming control: residence probability and $(D, T)$-stability,
{\it Automatica}, 28(3): 549-555, 1992.

\bibitem{Husher:1966}
H.~J.~Kushner, Finite time stochastic stability and the analysis of tracking systems,
{\it IEEE Trans. Automat. Contr.}, 11(2), 219-227, 1966.

\bibitem{Husher:1967}
H.~J.~Kushner, Stochastic Stability and Control, Academic Press, New York, 1967.




\bibitem {Liptser:1989}
R.~Sh.~Liptser and A.~N.~Shiryayev, Theorey of Martingales, Kluwer Academic, 1989.





\bibitem{Mao:2008}
X.~Mao, Stochastic Differential Equations and Applications, Second Edition, Woodhead, Cambridge, 2008.



\bibitem {Mao:2001}
X.~Mao, Some contributions to stochastic asymptotic stability and boundedness via multiple Lyapunov functions, {\it J. Math. Analy. Appl.}, 260(2), 325-340, 2001.


\bibitem{Meerkov:1988}
S.~M.~Meerkov and T.~Runolfsson, Residence time control,
     {\it IEEE Trans. Automat. Contr.}, 33(4), 323-332, 1988.


\bibitem{Bernt:2005}
B.~ {\O}ksendal, Stochastic Differential Equations: An Introduction with Applications, Sixth Edition, Springer-Verlag, New York, 2005.



\bibitem {Yin:2011}
J.~Yin, S.~Khoo, Z.~Man, and X.~Yu, Finite-time stability
and instability of stochastic nonlinear systems, {\it Automatica},
47(12): 2671-2677, 2011.
\bibitem{Yin:2014}
J.~Yin and S.~Khoo, Continuous finite-time state feedback stabilizers for
some nonlinear stochastic systems, {\it Int. J. Robust Nonlinear Control}, 25(11): 1581-1600, 2015.

\bibitem{Yu:2010}
X.~Yu and X.~Xie, Output feedback regulation of stochastic nonlinear
systems with stochastic iISS inverse dynamics, {\it IEEE
Trans. Automat. Contr.}, 55(2): 304-320, 2010.


\bibitem {Zak:2003}
S.~H.~ \.{Z}ak, Systems and Control, Oxford University Press, Oxford, 2003.

\end{thebibliography}
\end{document}